# Pellet Production Optimization using a Parallelized Progressive Hedging Algorithm

**Amin Aghalari**[1], **Badr Saleh Aladwan**[1], **Bruno Silva**[2], **Shaun Tanger**[3], **Mohammad Marufuzzaman**[1], **Veera Gnaneswar Gude**[4]

[1]Industrial & Systems Engineering, Mississippi State University, Starkville, MS, 39762
[2]Forest Management and Economics, Mississippi State University, Starkville, MS, 39762
[3]Coastal Research & Extension Center, Mississippi State University, Starkville, MS, 39762
[4]Civil and Environmental Engineering, Mississippi State University, Starkville, MS, 39762

**Abstract**

Renewable energy policies have driven the wood pellet market over the last decades worldwide. Among other factors, the return from this business depends largely on how well the producers manage the uncertainty associated with biomass yield and quality. This study develops a two-stage stochastic programming model that optimizes different critical decisions (e.g., harvesting, storage, transportation, quality inspection, and production decisions) of a biomass-to-pellet supply system under biomass yield and quality uncertainty. The goal is to economically produce pellets while accounting for the different pellet standards set forward by the U.S. and European markets. We propose two parallelization schemes to efficiently speed up the convergence of the overall decomposition method. We use Mississippi as a testing ground to visualize and validate the performance of the algorithms.

## 1  Introduction

The increasing commitment of government and private entities to reducing the carbon footprint have boosted the demand for wood pellet worldwide. The world pellet market is estimated to reach nearly 54 million tons by 2025, where approximately 70.4% of the total market demand is expected to be consumed in Europe alone [1]. These opportunities have leaned investors towards developing/improving the new/existing biomass-to-pellet supply chain, which can economically produce and transport pellets locally and for the overseas markets. Unlike other supply chains, the optimization of a biomass-to-pellet supply chain is extremely challenging. A number of factors drive towards letting the system optimization complicated, including but not limited to biomass quality variability, seasonality, market-specific pellet pro- duction requirements, and many others. Biomass supply chain has a rich literature. Most of the past studies tend to focus on minimizing the overall biomass supply chain costs (e.g., feedstock collection, inventory, production, and facility location decisions) under both deterministic (e.g., [2]) and stochastic (e.g., [3]) conditions. Even though the biomass-to-biofuel supply chain received considerable attention in the research community, the promising biomass-to-pellet supply chain is vastly unexplored. Until recently, a number of studies performed a techno-economic analysis to find an economical pathway to produce pellets from various agricultural sources [4, 5]. For instance, Agar [4] perform an economic comparison between the torrefied pellets and the conventional wood pellets. Hoefnagels et al. [5] perform an economic assessment to produce pellets from the underutilized wood resources available in the Southeast US. Based on the lo- cation and feedstock supply assumptions, the authors also estimate the optimal pellet plant size (55 to 315 Gg/year) for this region. A recent study (e.g., [6]) accounted the stochasticity associated with the biomass quality into the modeling process. Different from the study discussed above, our study, for the first time in the literature, aimed to focus on realistically modeling the key drivers into a biomass-to-pellet based sup- ply system network. The major contribution of this study is to first develop a mathematical model which explicitly captures the impact of biomass quality on a biomass-to-pellet supply network. The second contribution is the development of two parallelization schemes to efficiently speed up the convergence of the overall algorithm.

## 2 Problem Description and Model Formulation

This section presents a mathematical model that accounts for all the critical steps that significantly impact the biomass-to-pellet supply chain. Below is a summary of the sets, parameters, and decision variables of the optimization model.

**Sets:**
- $I$: set of feedstock suppliers, $i \in I$
- $J$: set of potential location for depots, $j \in J$
- $B, B_p$: set of biomass types, $b \in B$, and subset of possible feedstocks to produce pellet $p \in P$, respectively
- $P_A, P_E, P_b$: set of all pellet types, $p \in P$, subset of pellets to produce for the U.S. market, subset of pellets to produce for the Europe market, and subset of pellets that can be produced using biomass $b \in B$, respectively
- $C$: set of depot capacities, $c \in C$
- $T$: set of time periods, $t \in T$
- $R, R_p$: range of ash contents for pellets, and subset of ranges associated with each pellet type, respectively
- $\Omega$ : set of scenarios, $\omega \in \Omega$

**Parameters:**
- $\psi_{cj}$ : investment cost to open a depot
- $l_{bij\tau t\omega}$: unit transportation cost of flowing feedstock
- $h^1_{bi\tau t}, h^2_{bj\tau t}$: unit storage cost for feedstock type at supply site and depot, respectively
- $\sigma_{bit}$: unit cost of harvesting biomass
- $\zeta_{bjrt}$: unit cost inspection and segregation
- $\eta_{pjt}$: unit pellet production cost
- $\theta_{bjrr't}$: ash content reduction cost of unprocessed biomass
- $\pi_{pt}$: unit penalty cost for unsatisfied market demand of pellet
- $\bar{S}_{bit\omega}$ : availability of feedstock, and demand for pellets at a specific time
- $\alpha_{b\tau t\omega}$: dry matter loss (%) due to storing feedstock
- $\phi_{brp}$: conversion rate of unprocessed feedstock
- $\Gamma_{cj}$: pellet production capacity of specific size at a depot
- $\Pi^1_{bi}, \Pi^2_{bcj}$: storage capacity of feedstock in supply site and depot, respectively
- $\gamma_{bjt\omega}, \mu_{bjt\omega}$: inspected ash content (%) and the excess moisture content (%) in unprocessed biomass, respectively
- $I_{bjrt\omega}$: inspected ash content (%) and the excess moisture content (%) in unprocessed biomass, respectively
- $\phi^f_{b\tau t\omega}, \phi^v_{b\tau t\omega}$: dry matter loss (in %) during loading/unloading and transporting
- $\kappa_{bjrt}$: unit moisture content reduction cost for unprocessed biomass
- $\beta_{bjrt\omega}$: an indicator to denote if the preheating for the sampled biomass is required or not, i.e., $\beta_{bjrt\omega} = 1$, if preheating is required; 0 otherwise
- $\rho_\omega$: probability of scenarios

**Decision Variables**
- $Y_{cj}$: 1 if a depot of specific capacity is opened in potential location; 0 otherwise
- $S_{bit\omega}$: amount feedstock harvested in supply site at a time
- $X_{bij\tau t\omega}$: amount of transported biomass between two time periods between a supply site and a depot
- $H^1_{bi\tau t\omega}, H_{bj\tau t\omega}$: amount of feedstock stored at supply site and depot, respectively, between two time periods
- $P_{bj\tau t\omega}$: amount of unprocessed feedstock which are transported to depot and are ready for quality inspection
- $Z_{bjt\omega}$: total availability of unprocessed biomass that are ready for quality inspection at depot
- $Q_{bjrt\omega}$: amount of inspected unprocessed biomass at a depot in a time period with an specific ash content level
- $L_{brpjt\omega}$: amount of biomass with specific ash content level used to produce pellets
- $R_{bjrr't\omega}$: amount of unprocessed biomass needs to be adjusted concerning the ash level

- $D_{pjrt\omega}$: amount of pellets produced at depot with a specific ash content level
- $U_{ptw}$: shortage of pellets

With this, we are now ready to introduce formulation (BQP) below.

$$(BQP) \quad \text{Minimize} \sum_{c \in C} \sum_{j \in J} \psi_{cj} Y_{cj} \tag{1}$$

$$+ \sum_{\omega \in \Omega} \rho_\omega \sum_{t \in T} \left\{ \sum_{b \in B} \left( \sum_{i \in I} \sigma_{bit} S_{bit\omega} + \sum_{i \in I} \sum_{j \in J} \sum_{\tau \leq t} l_{bij\tau t\omega} X_{bij\tau t\omega} + \sum_{i \in I} \sum_{\tau \leq t} h^1_{bi\tau t} H^1_{bi\tau t\omega} \right. \right.$$

$$+ \sum_{j \in J} \sum_{\tau \leq t} h^2_{bj\tau t} H^2_{bj\tau t\omega} + \sum_{p \in P} \left( \sum_{j \in J} \eta_{pjt} D_{pjt\omega} + \pi_{pt} U_{pt\omega} \right)$$

$$\left. \left. + \sum_{j \in J} \sum_{r \in R} \left( \zeta_{bjrt} Q_{bjrt\omega} + \sum_{r' \leq r} (\theta_{bjrr't} + \kappa_{bjr't} \beta_{bjr't\omega} \mu_{bjt\omega}) R_{bjrr't\omega} \right) \right) \right\}$$

Subject to,

$$\sum_{c \in C} Y_{cj} \leq 1 \qquad \forall j \in J \tag{2}$$

$$S_{bit\omega} \leq \bar{S}_{bit\omega} \qquad \forall b \in B, i \in I, t \in T, \omega \in \Omega \tag{3}$$

$$S_{bit\omega} = H^1_{bitt\omega} + \sum_{j \in J} X_{bijtt\omega} \qquad \forall b \in B, i \in I, t \in T, \omega \in \Omega \tag{4}$$

$$(1 - \alpha_{b\tau,t-1,\omega}) H^1_{bi\tau,t-1,\omega} = H^1_{bi\tau t\omega} + \sum_{j \in J} X_{bij\tau t\omega} \qquad \forall b \in B, i \in I, (\tau, t) \in T | \tau \leq t - 1 \in T, \omega \in \Omega \tag{5}$$

$$\sum_{i \in I} X_{bijtt\omega} = H^2_{bjtt\omega} + P_{bjtt\omega} \qquad \forall b \in B, j \in J, t \in T, \omega \in \Omega \tag{6}$$

$$(1 - \alpha_{b\tau,t-1,\omega}) H^2_{bj\tau,t-1,\omega} + \sum_{i \in I} X_{bij\tau t\omega} = H^2_{bj\tau t\omega} + P_{bj\tau t\omega} \qquad \forall b \in B, j \in J, (\tau, t) \in T | \tau \leq t - 1 \in T, \omega \in \Omega \tag{7}$$

$$\sum_{\tau=1}^{t} H^1_{bi\tau t\omega} \leq \Pi^1_{bi} \qquad \forall b \in B, i \in I, t \in T, \omega \in \Omega \tag{8}$$

$$\sum_{\tau=1}^{t} H^2_{bj\tau t\omega} \leq \sum_{c \in C} \Pi^2_{bcj} Y_{cj} \qquad \forall b \in B, j \in J, t \in T, \omega \in \Omega \tag{9}$$

$$\sum_{\tau=1}^{t} P_{bj\tau t\omega} = Z_{bjt\omega} \qquad \forall b \in B, j \in J, t \in T, \omega \in \Omega \tag{10}$$

$$Z_{bjt\omega} = \sum_{r \in R} I_{bjrt\omega} Q_{bjrt\omega} \qquad \forall b \in B, j \in J, t \in T, \omega \in \Omega \tag{11}$$

$$\sum_{r \in R} Q_{bjrt\omega} = Z_{bjt\omega} \qquad \forall b \in B, j \in J, t \in T, \omega \in \Omega \tag{12}$$

$$\sum_{r' \leq r} R_{bjrr't\omega} = Q_{bjrt\omega} \qquad \forall b \in B, j \in J, r \in R, t \in T, \omega \in \Omega \tag{13}$$

$$\sum_{r' \leq r} R_{bjrr't\omega} = \sum_{p \in P_b} L_{br'pjt\omega} \qquad \forall b \in B, j \in J, r' \in R, t \in T, \omega \in \Omega \tag{14}$$

$$\sum_{b \in B_p} \phi_{br'p} L_{br'pjt\omega} = D_{pjr't\omega} \quad \forall p \in P, j \in J, r' \in R_{rp}\ t \in T, \omega \in \Omega \quad (15)$$

$$\sum_{p \in P} \sum_{r \in R_p} D_{pjr't\omega} = \sum_{c \in C} \Gamma_{cj} Y_{cj} \quad \forall j \in J, t \in T, \omega \in \Omega \quad (16)$$

$$\sum_{j \in J} \sum_{r \in R_p} D_{pjrt\omega} + U_{pt\omega} = d_{pt} \quad \forall p \in P, t \in T, \omega \in \Omega \quad (17)$$

$$Y_{cj} \in \{0,1\} \quad \forall c \in C, j \in J \quad (18)$$

$$X_{bij\tau t\omega}, H^1_{bi\tau t\omega}, H^2_{bj\tau t\omega}, P_{bj\tau t\omega} \geq 0 \quad \forall b \in B, i \in I, j \in J, (\tau, t) \in T | \tau \leq t - 1 \in T, \omega \in \Omega \quad (19)$$

$$S_{bit\omega}, Z_{bjt\omega}, Q_{bjrt\omega}, R_{bjrr't\omega}, D_{pjrt\omega}, U_{ptw} \geq 0 \quad \forall b \in B, p \in P, i \in I, j \in J, t \in T, (r, r') \in R | r' \leq r, \omega \in \Omega \quad (20)$$

The objective function (1) is the sum of the first-stage costs and the expected second-stage costs. Constraints (2) ensure that at most one depot of specific size is opened in a location. Constraints (3) restrict the harvesting quantity of biomass at supply sites in time to its availability. Constraints (4) and (5) are the flow balance constraints for biomass types at feedstock supply sites. Likewise, constraints (6) and (7) are the flow balance constraints for biomass types at depots. Constraints (8) and (9) set the biomass storage capacity at feedstock supply sites and depots, respectively. Constraints (10) determine the total biomass availability for quality inspection at depots. Constraints (11)-(14) classify the inspected biomass based on different ash levels. Constraints (15) ensure the conversion of raw biomass to pellets. Constraints (16) limit the pellet production capacity in depots. Constraints (17) ensure that the demand for pellets must be satisfied either via the depots or from a third-party supplier. Constraints (18) set the binary restrictions for the depot opening decisions. Finally, constraints (19) and (20) set the standard non-negativity restrictions.

## 3 Solution Approach

By setting $|B| = |P| = |T| = |R| = |\Omega| = 1$, model **(BQP)** can be reduced to a fixed charge network flow problem, which is already known to be an *NP*-hard problem [7]. As such, this study proposes a hybrid algorithm which combines the Sample Average Approximation (SAA) technique [8, 9] with an enhanced Progressive Hedging Algorithm (PHA) [10] to solve **(BQP)** in a reasonable timeframe. To further improve the computational efficiency of the proposed algorithm, two parallel schemes are developed. We employ two different parallelization schemes, referred to as Scheme 1 and 2, to efficiently solve the SAA replications. Both the enhancements exploit the multiprocessing capabilities of the local computers to efficiently solve **(BQP)**. **Scheme 1** assigns scores to each replication, E, based on the observed supply scenarios ($\bar{S}_{bitn}$) at the master node. The load and the assignment on the slaves (processors) are then made based on the scores. The scores $r^E$ are then sorted in descending order and assign the first p of them to *p* slaves (based on the available processors) in batches. Note that every time the master assigns a batch of *p* number of replications (based on the assigned scores) to the *p* slaves and the process continues until all the replications are solved. **Scheme 2** asynchronously assigns the SAA replications, *E*, to different available processors. The master is responsible for creating an initial pool for the SAA subproblems, where the subproblems are dynamically assigned to the available processors without following any specific order. The goal is to minimize the waiting time for the processors, which eventually speeds up the computation time and convergence of the overall SAA algorithm.

## 4 Computational Study and Managerial Insights

In this study, using Mississippi as a testing ground, a real-life case study is presented to visualize and validate the modeling results s. Finally, the computational performance of the proposed solution approaches is compared under varying test instances. All the algorithms are coded in Python 3.7 and executed on a desktop computer with Intel Core i7 3.60 GHz processor and 32.0 GB RAM. The optimization solver used is GUROBI 9.0.3.

### 4.1 Performance Evaluation of the Algorithms

This section presents the algorithms' computational performance in solving model **(BQP)**. In order to generate the test instances for performance evaluation of the algorithms, we vary sets $|I|, |J|$, and $|T|$. By doing so, 9 test instances with varying sizes are generated. To terminate the algorithms, the following termination criteria are used: (i) the optimality gap falls below a threshold value (e.g., $\epsilon = 1.0\%$); or (ii) the maximum time limit ($t^{max} = 1,800$) is reached ; or (iii) the maximum iteration limit ($r^{max} = 50$) is reached. The following notations are used to represent each particular variant of the proposed algorithms.
- **PHA**: Progressive Hedging Algorithm.
- **PHA+HR**: Enhanced Progressive Hedging Algorithm with application of Heuristics strategies.
- **PHA+HR+SB**: Enhanced Progressive Hedging Algorithm with application of both Heuristics strategies and Scenario Bundling techniques.
- **SAA**: Sample Average Approximation Algorithm.
- **Hybrid**: Hybrid decomposition algorithm combining Sample Average Approximation and Enhanced Progressive Hedging Algorithm (PHA+HR+SB).
- **Hybrid+PL1**: Parallelization scheme I is applied over hybrid algorithm **Hybrid**.
- **Hybrid+PL2**: Parallelization scheme II is applied over hybrid algorithm **Hybrid**.

The results reported in Table 1 signifies that the computational performance of the basic **PHA** can be further improved by incorporating different accelerating techniques. Results in Table 2 further demonstrate that incorporating **PHA+HR+SB** under a **SAA** framework, i.e., the **Hybrid** algorithm, significantly reduces the average solution time for solving different instances of **(BQP)**. The results reported in Table 3 indicate that parallelizing the replications of the **SAA** algorithm using *scheme 2* saves approximately 15.6% solution time over algorithm **Hybrid+PL1**, while both provide competitive quality solutions.

Table 1: Performance of the **PHA** algorithm

|  | Gurobi | | PHA | | | PHA+HR | | | PHA+HR+SB | | |
| --- | --- | --- | --- | --- | --- | --- | --- | --- | --- | --- | --- |
|  | $t(sec)$ | $\epsilon(\%)$ | $t(sec)$ | $\epsilon(\%)$ | $r$ | $t(sec)$ | $\epsilon(\%)$ | $r$ | $t(sec)$ | $\epsilon(\%)$ | $r$ |
| Average over 9 instances | 10,610 | 11.23 | 11,107 | 10.23 | 10 | 10,625 | 5.23 | 7 | 9,548 | 1.34 | 5 |

Table 2: Performance of **PHA+HR+SB**, **SAA**, and **Hybrid** algorithm

|  | PHA+HR+SB | | | SAA | | | Hybrid | | |
| --- | --- | --- | --- | --- | --- | --- | --- | --- | --- |
|  | $t(sec)$ | $\epsilon(\%)$ | $r$ | $t(sec)$ | $\epsilon(\%)$ | $r$ | $t(sec)$ | $\epsilon(\%)$ | $r$ |
| Average over 9 instances | 9,548 | 1.34 | 5 | 9,181 | 3.1 | 2 | 7,653 | 1.66 | 2 |

Table 3: Performance of **Hybrid** algorithm under different parallelization schemes

|  | Hybrid | | | Hybrid+PL1 | | | Hybrid+PL2 | | |
| --- | --- | --- | --- | --- | --- | --- | --- | --- | --- |
|  | $t(sec)$ | $\epsilon(\%)$ | $r$ | $t(sec)$ | $\epsilon(\%)$ | $r$ | $t(sec)$ | $\epsilon(\%)$ | $r$ |
| Average over 9 instances | 7,653 | 1.66 | 2 | 6,146 | 0.54 | 2 | 5,186 | 0.48 | 2 |

### 4.2 Case Study

We now present a real-world case study to understand the impact of biomass quality variability (e.g., ash or moisture content) on pellet production and network design decisions. To serve this purpose, we create five different biomass quality scenarios, as illustrated in Table 4.

Table 4: Description of biomass quality scenarios

| Quality Scenario | Description |
| --- | --- |
| Base | Assuming the base *ash* and *moisture* qualities |
| Low ash | Base *moisture* but *ash* quality dropped by 30% from the base case |
| Good ash | Base *moisture* but *ash* quality improved by 30% from the base case |
| Low moisture | Base *ash* but *moisture* quality dropped by 30% from the base case |
| Good moisture | Base *ash* but *moisture* quality improved by 30% from the base case |

Results in Table 5 delineate that pellet depot location decisions ($|J|$) are sensitive to the biomass quality parameters. We observe that the base quality case decides to open 34 pellet facilities within our test region. However, as the biomass quality improves, model **(BQP)** decides to open 36 and 35 depots under good ash and moisture quality scenarios, respectively. Overall, we observe that high-quality biomass promotes the growth of pelleting businesses; however, decision-makers should carefully consider the local biomass quality parameters to make appropriate site selection decisions.

Table 5: Number of depots selected under different biomass quality scenarios

| Base | Ash Quality | | Moisture Quality | |
| --- | --- | --- | --- | --- |
|  | Good | Low | Good | Low |
| 34 | 36 | 31 | 35 | 33 |

Figure 1 illustrates how different biomass quality parameters impact the biomass storage decisions (at the pelleting facilities). The results reveal that storage is considered a favorable option when the biomass quality is high, while the converse is true when the quality of the biomass is low. For instance, it is observed that as the biomass base ash quality shifts to low- and high-quality scenarios, the overall biomass storage is changed by -19.8% and 20.1%, respectively, from the base case. A similar trend is observed for the moisture content for which the decisions are changed by -45.6% and 41.3%, respectively, from the base moisture quality scenario. Finally, we observe that the storage decisions are more sensitive to the moisture content than the ash content.

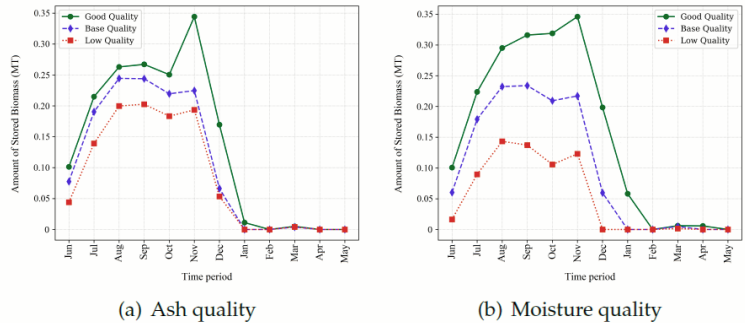

Figure 1: Biomass storage decisions under different quality levels

## 5 Conclusions and Future Research Directions

This paper develops a two-stage stochastic mixed-integer linear programming model to capture several realistic features biomass-to-pellet supply chain. We then develop an innovative solution approach, enhanced by the parallel computing concepts, to efficiently solve the proposed optimization model. We assume that the system is reliable and will never fail. It might be interesting to investigate the sudden surge of woody biomass due to a natural catastrophe (e.g., hurricane) or human-induced catastrophe [11, 12, 13] and its association with the pellet production decisions.